\documentclass[9pt,twocolumn]{IEEEtran}
\usepackage{amsmath} 
\usepackage{amssymb}  
\usepackage{amsfonts}
\usepackage{graphicx}
\usepackage{subfigure}
\usepackage[usenames]{color}

\newtheorem{remark}{Remark}

\newtheorem{dfn}{Definition}
\newtheorem{thm}{Theorem}
\newtheorem{lem}{Lemma}
\newtheorem{prop}{Proposition}

\newenvironment{pf}{\smallbreak\noindent{\it Proof. }}{\hfill$\Box$\smallbreak}
\newenvironment{pf*}[1]{\smallbreak\noindent{\it #1}}{\hfill$\Box$\smallbreak}

\DeclareMathOperator{\trace}{tr}

\newcommand{\RR}{\mathbb{R}}
\newcommand{\CC}{\mathbb{C}}
\newcommand{\ZZ}{\mathbb{Z}}
\newcommand{\NN}{\ZZ_+}

\renewcommand{\d}{\delta}

\renewcommand{\t}{\theta}
\newcommand{\e}{\epsilon}
\newcommand{\D}{\Delta}

\newcommand{\sA}{\mathcal{A}}

\newcommand{\SE}{J^{SE}}

\newcommand{\dx}{{n_x}}

\newcommand{\dy}{{n_y}}

\newcommand{\dw}{{n_w}}
\newcommand{\dz}{{n_z}}

\newcommand{\X}{\RR^{\dx}}

\newcommand{\W}{\RR^{\dw}}

\newcommand{\cvgip}{\stackrel{\textrm{P}}{\rightarrow}}

\begin{document}

\title{\LARGE \bf
Stable Nonlinear Identification From Noisy Repeated Experiments via Convex Optimization
}
\author{Mark M. Tobenkin$^{\dagger,+}$,
   Ian R. Manchester$^{\star}$,
  and Alexandre Megretski$^{\dagger}$
\thanks{*Supported by National Science Foundation Grant
  No. 0835947.}
\thanks{$^{+}$ Corresponding author.}
\thanks{$^{\star}$ School of Aerospace, Mechanical, and Mechatronic
  Engineering, University of Sydney, Sydney, NSW, 2006. }
\thanks{$^{\dagger}$ Department of Electrical Engineering and Computer Science at the Massachusetts Institute of Technology, Cambridge, MA 02149.}
}

\maketitle

\begin{abstract}
This paper introduces new techniques for using convex optimization to fit input-output data to a class of stable nonlinear dynamical models. We present an algorithm that guarantees consistent estimates of models in this class when a small set of repeated experiments with suitably independent measurement noise is available. Stability of the estimated models is guaranteed without any assumptions on the input-output data. We first present a convex optimization scheme for identifying stable state-space models from empirical moments. Next, we provide a method for using repeated experiments to remove the effect of noise on these moment and model estimates. The technique is demonstrated on a simple simulated example.

\end{abstract}


\section{INTRODUCTION}
Building nonlinear dynamical models capable of accurate long term prediction is a common goal in system identification. However, for most model structures multi-step prediction errors have a complex nonlinear dependence on the model parameters.  
Furthermore, 
assuring stability of algorithmically generated nonlinear models
is a substantial challenge. 
In many practical situations, where data-sets are limited or under-modeling is present, widely used ``one-step'' prediction error minimization techniques can render models that are unstable or have poor multi-step predictions.  This work presents a convex optimization 
method for approximating the input-output response of a nonlinear dynamical system via state-space models with stability guarantees. This paper extends recent work in \cite{Tobenkin10}, \cite{Bond10} and \cite{Megretski08} by providing a family of consistent estimators for a class of stable nonlinear models when a small set of repeated experiments 
is available.
We examine the problem of embedding an input-output identification task inside 
a state-space modeling framework.  We inherit from the methods 
of \cite{Tobenkin10}, \cite{Bond10}, \cite{Megretski08}
an unqualified guarantee of model stability and a cost function that is a convex upper bound on the ``simulation error'' associated with these models.  However, the estimators from
\cite{Tobenkin10}, \cite{Bond10}, \cite{Megretski08}
are generally not consistent, and for systems that are nearly marginally stable the biasing effect of measurement noise can be quite severe.  Furthermore, the complexity of these methods grows undesirably with the number of data points.

We present a modification of algorithms from \cite{Tobenkin10} that mitigates these two difficulties.  In particular, a 
technique that utilizes the problem data through empirical moments
only is used. 
As a result, the complexity of the method generally grows linearly with data-set size.  We also provide a method  for asymptotically removing the effects of measurement noise on these empirical moments when a small set of repeated experiments are available, utilizing an idea which is superficially similar to 
instrumental variable methods \cite{Ljung99}.  We that demonstrate that this technique, a nonlinear extension of \cite{Manchester12}, recovers consistency when the data is generated by a system within a specific class of models.  


\subsection{Previous Work}
The use of maximum likelihood and one-step prediction error methods is frequently motivated by the consistency and  asymptotic efficiency of the resulting estimators \cite{Ljung99}.  In the face of limited data or significant under-modeling, these techniques often render models that are unstable or make poor multi-step ahead predictions \cite{Farina10}.  Direct minimization of longer term prediction errors have appeared in several forms, including the output-error method for input-output system identification, \cite{Soderstrom82}, notions of ``best'' approximation, \cite{Paduart10}, and simulation error minimization, \cite{Bonin10},\cite{Farina10}.  These methods require optimization of a  non-convex functional for all but the simplest model structures (e.g. finite impulse response and Volterra type models) and can suffer from local minima \cite{Soderstrom82}. Appealing theoretical properties of these methods (e.g. efficiency and unbiasedness) are often predicated on finding global minima of generically hard nonlinear programming problems.

Several results are available for linear time invariant (LTI) system identification using least squares that provide stability guarantees even in the face of under-modeling (e.g.  \cite{Regalia94},\cite{Regalia95}, \cite{Tugnait98}).  It is worth noting that these stability guarantees apply only as the number of available data points tends to infinity and 
requires an assumption that the data is generated by a (potentially under-modeled) stationary LTI process.  Several modified subspace techniques have also been presented to address the issue of model stability.  In \cite{Gestel01} regularization is used to ensure model stability.  In \cite{Lacy02} and \cite{Hoagg04} a joint search over Lyapunov function and system dynamics using convex optimization was used to ensure model stability.  The LTI-specific
method employed by 
\cite{Lacy02} and \cite{Hoagg04}
is closely related to the technique by which this 
paper addresses stability.

Several convex relaxation techniques have recently been employed by the Set Membership (SM) identification community to address fixed order identification of LTI systems (\cite{Cerone11}, \cite{Feng11}).  In \cite{Cerone11} outer approximations of the set of parameters consistent with bounded noise and stability assumptions are computed.  In \cite{Feng11} a convex relaxation approach is suggested for optimization of arbitrary polynomial objectives over the set of LTI models consistent with a given data-set and a set of stability and bounded noise assumptions. A similar approach is taken for identifying Linear Parameter Varying systems in \cite{Cerone12}. By contrast, 
in this work we examine a ``convex restriction'' 
approach where inner approximations of the set of stable models are used to guarantee stability and convex upper bounds on the cost function of interest are used as a surrogate objective.

\subsection{Outline}
The paper proceeds as follows. Section~\ref{sec:preliminaries} presents the notation, problem setup, and a bias elimination strategy employed in this work.  Next, Section~\ref{sec:convex} provides a convex parameterization of stable state-space models and a convex upper bound for simulation error.  This parameterization and objective are then combined with the bias elimination strategy in Section~\ref{sec:noise}, wherein a system identification algorithm based on semidefinite programming is given along with asymptotic analysis of the method.  Finally, a comparison of the proposed algorithm to two alternative least-squares based methods is provided in Section~\ref{sec:example}.

\section{Preliminaries}
\label{sec:preliminaries}
In this section we introduce basic notation, and present the problem
setup to be addressed in the paper.

\subsection{Notation}
\label{sec:notation}

$\CC^{k\times n}$ stands for the set of all $k$-by-$n$ 
complex matrices, with $\CC^n$ being a shorthand for 
$\CC^{n\times 1}$.
$\RR^{k\times n}$ and $\RR^n$ are the subsets of real matrices from
$\CC^{k\times n}$ and $\CC^n$ respectively. 
$\NN^n$ is the subset of $\RR^n$ whose elements are
non-negative integers. 

We use some notation from MATLAB, where
$A'$, $[A,B]$, and $[A;B]$ denote, respectively,
 Hermitian conjugation,
horizontal concatenation, and vertical concatenation of matrices.
For $R \in \CC^{k\times n}$ we denote by 
$[R]_{a,b}$
the scalar element in the $a$-th row and $b$-th column of $R$,
with the shorthand $[v]_d=v_{d,1}$ used for 
$v \in \CC^n$. In addition, for $v \in \CC^n$ and
$\alpha \in \NN^n$, 
\[ v^\alpha := \prod_{d=1}^n [v]_d^{[\alpha]_d}\]
is the monomial function of $v$ with {\sl vector degree} $\alpha$,
and {\sl scalar degree} $\|\alpha\|_1$, where, 
for $w\in\mathbb C^n$, $\|w\|_1 := \sum_{i=1}^n |[w]_i|$ is
the $\ell_1$
norm of $w$.  
For Hermitian matrices $A,B \in \CC^{n \times n}$
(i.e. such that $A=A'$ and $B=B'$),
$A \geq B$ (or $A > B$) means that $A-B$ is
positive semidefinite
(respectively, positive definite).
For $R = R' \in\CC^{n\times n}$ and
$v\in\CC^n$ we use the shorthand $|v|_R^2=v'Rv$. Moreover, when
$R\ge0$, we also write $|v|_R := \sqrt{v'Rv}$.  When $W$ is a set,
$\ell_T(W)$ denotes the set of all functions 
$w:~\{0,1,\dots,T\}\to W$. Naturally, the elements of  $\ell(W)$
are finite length {\sl sequences} of elements from $W$. 
The notation $\cvgip$ refers to
convergence in probability. 
\subsection{Problem Setup}
\label{sec:prb}
We define a {\sl data set with $N$ experiments of length $T$,
$n_w$-dimensional input, and
$n_x$-dimensional state} as a collection
$(\tilde w,\tilde x_1,\dots,\tilde x_N)$ of
 sequences
$\tilde w\in\ell_T(\mathbb R^{n_w})$, 
$\tilde x_i\in\ell_T(\mathbb R^{n_x})$. 
$\mathcal D(n_x,n_w,N,T)$ stands for the set of all data sets
of given dimensions, number of experiments, and signal length.
Accordingly, 
$\mathcal D(n_x,n_w,N)=\cup_{T=0}^\infty
\mathcal D(n_x,n_w,N,T)$ stands for the set of all
data sets with unspecified signal length.
\begin{figure}
  \centering
  \includegraphics[width=0.80\columnwidth]{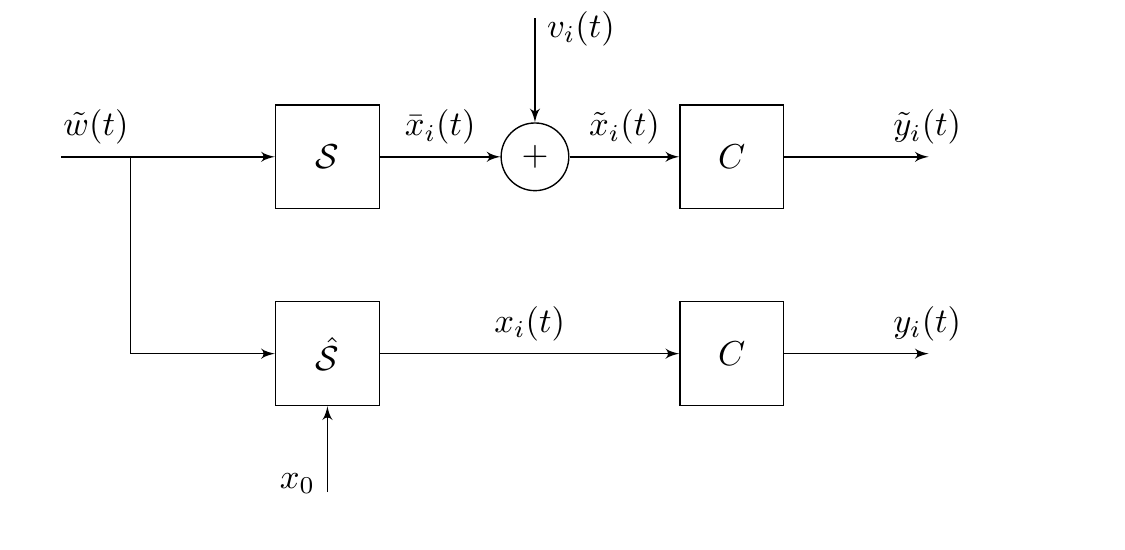}
  \caption{The experimental setup considered in this work. \label{fig:experiments}}
\end{figure}

In applications, each $\tilde x_i(t)$ is the result of
feeding the same input $\tilde w(t)$ into a system, ${\cal S}$, and
measuring the sum $\tilde x_i(t) = \bar x_i(t) + v_i(t)$, where $\bar
x_i(t)$ is the ``true system response'' and $v_i(t)$ is corrupting
measurement noise.  Additionally, in order to set a measure of 
quality for model predictions, we define an output signal $\tilde y_i(t)$ is defined
by $\tilde y_i(t) = C\tilde x_i(t)$, for some fixed matrix $C \in
\RR^{\dy\times\dx}$.  
Informally, the identification objective will be to accurately
predict the input-output behavior of this a system with $\tilde w(t)$
taken as input and $\tilde y(t)$ taken as an output (alternatively $C$
can be seen as weighting the importance certain components of $\tilde x_i(t)$).
This experimental setup is depicted in Figure~\ref{fig:experiments}.

The underlying assumption is that the collection of signals $\bar
x_i(t)$ constitute a reasonable state, or reduced state, for a
state-space model approximating the system behavior.  As an example, when identifying a
SISO system with input $u = u(t)$ and output $y = y(t)$, one can
imagine feeding in $N\cdot D$ samples of a $D$-periodic input $\tilde u(t)$
and measuring $\tilde y(t) = \bar y(t) + v_y(t)$, where $\bar y(t)$ is
the true system response and $v_y(t)$ is measurement noise.  In this
case, one could use the above setup with $\dx = \dw = n < D$ by taking
\[
\tilde w(t) = \begin{bmatrix} \tilde u(t+n) \\ \vdots \\ \tilde u(t+1) \end{bmatrix},
\quad
\tilde x_{i+1} =  \begin{bmatrix} \tilde y(t+n+iD) \\ \vdots \\ \tilde y(t+1+iD) \end{bmatrix},
\]
\[
v_{i+1} =  \begin{bmatrix} v_y(t+n+iD) \\ \vdots \\ v_y(t+1+iD) \end{bmatrix},
\]
for $i \in \{0,\ldots, N-1\}$ and $t \in \{0,\ldots,D-n\}$.  Here, the
matrix $C \in \RR^{1 \times n}$ might be $C = \begin{bmatrix} 1 & 0 &
  \ldots & 0 \end{bmatrix}$.

\subsection{State Space Models}
In general, a nonlinear state space model (time
invariant, in discrete time) with 
$n_w$-dimensional input and
$n_x$-dimensional state is specified by a function 
$a:~\mathbb R^{n_x}\times\mathbb R^{n_w}\to\mathbb R^{n_x}$,
which in turn defines the {\sl input-output function}
$G_a:~\mathbb R^{n_x}\times\ell(\mathbb R^{n_w})\to
\ell(\mathbb R^{n_w})$ mapping {\sl initial state} 
$x_0\in\mathbb R^{n_x}$ and {\sl input sequence} 
$w\in\ell(\mathbb R^{n_w})$ to the {\sl output sequence}
$x\in\ell(\mathbb R^{n_x})$ according to
\begin{equation}\label{eq:dynamics} x(t)=a(x(t-1),w(t)),\qquad
  x(0)=x_0.\end{equation}
For $w: \ZZ_+ \mapsto \RR^{n_w}$, we define $x =
G_a(x_0,w)$ to be the sequence similarly
defined by this recurrence.
Let $x=G_a(x_0,w)$, $\hat x=G_a(\hat x_0,w)$ be two responses of
system \eqref{eq:dynamics} to the same input $w: \ZZ_+ \mapsto \RR^{n_w}$ and different
initial conditions $x_0,\hat x_0$.
We call system \eqref{eq:dynamics} 
{\sl $\ell^2$-incrementally stable}
when $x-\hat x$ is square summable for all $x_0,\hat x_0,w$.  The
system \eqref{eq:dynamics} is {\sl incrementally exponentially stable} if
there exist constants $c > 0$ and $\rho \in (0,1)$, independent of
$x_0,\hat x_0$ and $w$, such that
$|x(t)-\hat x(t)| \leq c \rho^t |x_0 - \hat x_0|$ for all $x_0,\hat
x_0, w$ and $ t \geq 0$.

This paper deals with subsets of 
state space models \eqref{eq:dynamics} which have more specific
finite dimensional structure. 
For positive integers $n_x,n_w,n_\theta$ let 
$\Theta,\Phi,\Psi$ be a non-empty set  
$\Theta\subset\mathbb R^{n_\theta}$ and two sequences
$\Phi=\{\phi_i\}_{i=1}^{n_\theta}$,
 $\Psi=\{\psi_i\}_{i=1}^{n_\theta}$ of real analytical functions
$\phi_i:~\mathbb R^{n_x}\times\mathbb R^{n_w}\to
\mathbb R^{n_x}$, $\psi_i:~\mathbb R^{n_x}
\to\mathbb R^{n_x}$. We say that the 3-tuple $(\Theta,\Phi,\Psi)$
is a {\sl stable projective parameterization
with $n_w$ inputs, $n_x$ states, and $n_\theta$ parameters}
when,
for all $\theta\in\Theta$,
\begin{itemize}
\item the function $e_\theta:~\mathbb R^{n_x}
\to\mathbb R^{n_x}$ defined by
$ e_\t(x) = \sum_{i=1}^{n_\theta}[\theta]_i\psi_i(x)$
is a bijection;
\item the state space model
\eqref{eq:dynamics} with
 $a=a_\theta$ defined by
\begin{equation}
  a_\theta(x,w)=e_\t^{-1}(f_\t(x,w))
\end{equation}
\begin{equation}
f_\t(x,w) = \sum_{i=1}^{n_\theta}[\theta]_i\phi_i(x,w), \quad 
e_\t(\xi) = \sum_{i=1}^{n_\theta}[\theta]_i\psi_i(\xi),
\label{eq:implicit_dynamics}
\end{equation}
is $\ell^2$-incrementally stable.
\end{itemize}
Once a stable projective parameterization $(\Theta,\Phi,\Psi)$
is selected, a stable  state space model can be defined
by specifying a vector parameter $\theta\in\Theta$.

Recent discussion and applications of incremental stability and the
related notions of contractive and convergent systems can be found in
\cite{Fromion99}, \cite{Lohmiller98}, \cite{Angeli02},
\cite{Sontag10}, and \cite{Besselink12}.  This property is related to
familiar ``fading memory'' conditions employed in other identification
and system approximation papers (e.g. \cite{Boyd85},
\cite{Partington96}, \cite{Pillonetto11}), though we note that the
condition as defined above does not restrict the behavior of solutions
with different input sequences.

In practice, both the particular construction of the signals $\tilde w(t)$ and $\tilde
x_i(t)$ from measureable quantities, and the selection of the sequences of functions
$\Psi$ and $\Phi$ defining $(e_\t,f_\t)$ should be
guided by  a model selection criteria such as cross-validation
\cite{Ljung99}.  We consider both of these selections fixed for the
remainder of the paper.

\subsection{Empirical Moments}
For given positive integers $n_x,n_w,N$ 
let $n_z=2n_x+n_w$ be the dimension of the vectors 
\begin{equation}\label{am2}\tilde z_i(t)=
\begin{bmatrix}\tilde x_i(t)\\ \tilde x_i(t-1)\\ \tilde w(t)
\end{bmatrix}\qquad(t\in\{1,\dots,T\},~i\in\{1,\dots,N\})
\end{equation}
defined by the data set $\Xi = (\tilde w,\tilde x_1,\dots,\tilde x_N)
\in {\cal D}(\dx,\dw,N,T)$.
For $\alpha\in\mathbb Z_+^{n_z}$ such that $\|\alpha\|_1\le N$
and for
$z_1,\dots,z_N\in\mathbb R^{n_z}$ define
\[  p_\alpha(z_1,\dots,z_N)=
\prod_{i=1}^{\|\alpha\|_1} [ z_i]_{\beta_\alpha(i)},
\]
where
\[
\beta(i) = \min \left\{ d \in\{1,\ldots,N\}\; \left | \;
  \sum_{j=1}^{d}[\alpha]_j \geq i\right .\right\}.
\]
By construction, $z^\alpha=p_\alpha(z,\dots,z)$, so that
one can view $p_\alpha$ as a multi-linear function which generates
the monomial $z^\alpha$ when evaluated on the multi-diagonal
(note that such multi-linear functions are not uniquely defined by
$\alpha$).
For a given data set $\Xi=(\tilde w,\tilde x_1,\dots,\tilde x_N)
$ and $\alpha\in\mathbb Z_+^{n_z}$
define the {\sl linearized empirical moment} 
$\tilde \mu_{\alpha}(\Xi)$ by
\begin{equation}
\tilde \mu_{\alpha}=
  \tilde \mu_{\alpha}(\Xi) = \frac{1}{T} \sum_{t=1}^{T}
p_\alpha(\tilde z_1(t),\dots,\tilde z_N(t)).
  \label{eq:tildemu}
\end{equation}
Since it is sometimes convenient to emphasize 
$\tilde \mu_{\alpha}(\Xi)$ as a function of variable $\alpha$ with
a fixed $\Xi$, we will also use the equivalent notation
$\tilde \mu_{\alpha}(\Xi)=\hat\mu_\Xi(\alpha)$.
According to this notation, for a given data set $\Xi$ with $N$
experiments, $n_x$ states, and $n_w$ inputs, $\hat\mu_\Xi$ is a
real-valued function defined on the set of elements $\alpha\in\mathbb
Z_+^{2n_x+n_w}$ such that $\|\alpha\|_1\le N$.

Informally speaking, linearized 
empirical moments represent an attempt at
``de-noising'' the data contained in the vectors $\tilde z_i(t)$,
as defined by \eqref{am2} in the case
when $\tilde x_i(t)= \bar x(t)+v_i(t)$ for $t\in\{0,1,\dots,T\}$, 
where $ \bar x= \bar x(t)$, the ``true system response'',
 does not depend on the
experiment number $i$, 
and the noise variables $v_i(t)$ are suitably independent of
$\bar x$ and of each other, to produce good estimates 
$\tilde \mu_{\alpha}(\Xi)$ of the standard empirical moments
\begin{equation}  \mu_{\alpha}( \bar x,\tilde w)=\frac1T\sum_{t=1}^T
 \bar z(t)^\alpha,\qquad\bar z(t)=
\begin{bmatrix} \bar x(t)\\  \bar x(t-1)\\ \tilde w(t)
\end{bmatrix}.
\end{equation}
This  approach is inspired by instrumental variable (IV) techniques,
\cite{Soderstrom83}, with repeated experiments playing a role
comparable to a specific choice of instruments.  Rather than
asymptotically approximate a least squares parameter estimate, as in
IV methods, this work focuses on asymptotically minimizing an
alternative convex loss function that depends only on empirical moments.  To have a meaningful
convergence of the linearized empirical moments we require both the
aforementioned independence of the noise sequences, to be made more
precise shortly, and that the
true system responses, $\bar x_i(t)$, tend to one another despite their
differing initial conditions.  

\subsection{Persistence of Excitation}
The following  notion of persistence of excitation will be
used in our consistency analysis.

\begin{dfn}\label{def:persist}
  Fix two signals  $w: \NN \to \W$ and $x: \NN \to \X$, and let $w^{(T)}$
  and $x^{(T)}$ be the restriction of these signals to $\{0,\ldots,T\}$.
  For a given function $a: \X \times \W \to \X$, we say a pair of
  signals $(w,x)$ is {\it persistently exciting for $a$} if
  there exists a positive measure $\pi$ on $\X\times \W$ such that
  $\pi $ is supported on an open set, and for every finite subset $\aleph =
  \{\alpha_j\}_{j=1}^{|\aleph|}$ of $\NN^\dz$
$$
\liminf_T \lambda_{\textrm{min}}(M_T^{\aleph} -
  M_{\pi}^{\aleph}) \geq 0,
$$
where $M_T^\aleph, M_{\pi}^{\aleph} \in
\RR^{|\aleph|\times|\aleph|}$ are defined by:
\[
[M_T^{\aleph}]_{i,j} =  \mu_{\alpha_i+\alpha_j}(w^{(T)},x^{(T)}),
\]
\[
[M_\pi^{\aleph}]_{i,j} =  \int [a(x,w) ; x ; w]^{\alpha_i+\alpha_j} d\pi(x,w).
\]
\end{dfn}

Informally, this non-standard notion of  persistence of excitation will be employed
to establish a connection between 
$$\frac{1}T \sum_{t=1}^T| e_\t(a(x(t),w(t))) - f_\t(x(t),w(t))|^2$$
vanishing as $T \rightarrow \infty$ and $a$ being equivalent to $a_\t
= e_\t^{-1}\circ f_\t$.  The use of a projective representation,
i.e. $a_\t$ being implicitly defined, renders several complications to
standard consistency arguments based on strong convexity (for example,
the $e_\t$ and $f_\t$ that define $a_\t$ can be non-unique).  The above
notion of persistence will be used to circumvent these difficulties.



\subsection{Data-Matching Error}
We examine the following loss function for identifying models.
\begin{dfn}
 The {\it T-step simulation error}, $\SE_T$, is a function of an $a:
 \X \times \W \to \W$, an initial condition vector
 $x_0 \in \X$, and two signals $\tilde w \in \ell_T(\W),$ and $ \tilde
 x \in
 \ell_T(\X)$, defined by
 \begin{equation}
    \SE_T(a,x_0,\tilde x, \tilde w) = \frac{1}{T}\sum_{t=0}^{T-1} |C(\tilde x(t) - x(t))|^2,
  \end{equation}
  where $x = G_{a}(x_0,\tilde w)$.
\end{dfn}

\subsection{Data Generation Mechanism}
\label{sec:mechanism}
Two data generation mechanisms, defined by  considering data sequences
as stochastic processes, will be analyzed in this work.  These
mechanisms  consider signals defined on
an infinite horizon, i.e.  $\tilde w:
\NN \to \W$ and $\tilde x_i: \NN \to \W$, for $i \in\{1,\ldots,N\}$.  We express $\tilde x_i$ as
the sum of two signals $\bar x_i$ and $v_i$, again representing the
true system response and measurement noise respectively.  Let $\tilde
x_i^{(T)}$,
and $\tilde w^{(T)}$ be the restrictions to $\{0,\ldots,T\}$ of $\tilde
x_i$, and $\tilde w$ respectively.  Then we define the data
set $\Xi_T \in {\cal D}(\dx,\dw,N,T)$ by
$$\Xi_T = (\tilde w^{(T)},\tilde x_1^{(T)},\ldots,\tilde x_N^{(T)}).$$

The following assumptions define the first data generation mechanism.
\begin{enumerate}
\item[(A1)] The signal $\tilde w(t)$ is a stochastic
  process for $t \in \ZZ$, which is uniformly bounded in $t$.
\item[(A2)] The signals $v_i(t)$ are i.i.d. zero mean bounded
  stochastic processes independent of one another, $\tilde w(t)$ and
  each $\bar x_i(t)$.
\item[(A3)] The signals $\bar x_i(t)$ are stochastic processes which
  are uniformly bounded in $i$ and $t$.  There exist constants $c > 0$ and $\rho \in
  (0,1)$ such that.
\begin{equation}
  |\bar x_i(t)-\bar x_j(t)| \leq c \rho^t, \qquad \forall i,j \in
  \{1,\ldots,D\}, t \in \NN,\label{eq:yconv}
\end{equation}
almost surely;
\end{enumerate}

An alternative, less general, data-generation mechanism in given by
the assumptions (A1), (A2), and the following.
\begin{enumerate}
\item[(A4)] There exists a function $a_0: \X \times
  \W \to \X$ such that \eqref{eq:dynamics} with $a = a_0$ defines a
  BIBO and  incrementally
  exponentially stable system and $\bar x_i = G_{a_0}(\bar
  x_{i0},\tilde w)$ for some unknown $\bar x_{i0} \in \X$, for $i \in \{1,\ldots,N\}$. The pairs of signals
  $(\tilde w, \bar x_i)$ are persistently exciting with respect to $a_0$,
  as in Definition~\ref{def:persist}, with probability one.
\end{enumerate}
The appendix contains practical conditions on $\tilde w(t)$ and $a_0$ that ensure (A4) holds.
It is immediate that assumptions (A1) and (A4) together imply (A3).

\subsection{Identification Objective}
In this paper, we view system identification algorithms as
4-tuples $(\mathcal A,\Theta,\Phi,\Psi)$, where 
$(\Theta,\Phi,\Psi)$ is a  stable projective parameterization
with $n_w$ inputs, $n_x$ states, and $n_\theta$ parameters,
and $\mathcal A$ is a function 
$\mathcal A:~\mathcal D(n_x,n_w,N)\to\Theta$ mapping
data sets to parameter vectors from $\Theta$.

Specifically, we are interested in generating efficient
 {\sl moments-based}
system identification algorithms $(\mathcal A,\Theta,\Phi,\Psi)$,
i.e. those
for which the function $\mathcal A:~\mathcal D(n_x,n_w,N)\to\Theta$
has the form $\mathcal A(\Xi)=\mathbf A(\hat\mu_\Xi)$,
which means that the resulting identified 
model is a function of the linearized
empirical moments $\hat\mu_\Xi(\alpha)$ with 
$\alpha\in\mathbb Z_+^{2n_x+n_w}$ satisfying $\|\alpha\|_1\le N$.

The main contribution of this paper is the construction of moments-based
system identification algorithms $(\mathcal A,\Theta,\Phi,\Psi)$ and
sets $\Theta_0 \subset \Theta$
with the following properties:
\begin{itemize}
\item[(a)] the set $a_{\Theta_0}=\{a_\theta:~\theta\in\Theta_0\}$
of models \eqref{eq:dynamics} generated by $\Theta_0$
is sufficiently broad,
in the sense that
every stable linear state space model $a(x,w)=Ax+Bw$ 
is in $a_{\Theta_0}$, and some non-linear functions are contained in
$a_{\Theta_0}$ as well;
\item[(b)] when a sequence of data sets $\{\Xi_T\}_{T=1}^\infty$ is
  generated by signals $(\tilde w, \tilde x_1, \ldots, \tilde x_N)$
  satisfying assumptions (A1),(A2) and (A3),  then $\theta_T = \sA(\Xi_T)$ asymptotically (with respect to $T$)
  minimizes an upper bound for
  \begin{equation}
    \frac{1}{N}\sum_{i=1}^N \SE_T(a_\t,\tilde x_0,\bar x_i^{(T)}, \tilde w^{(T)}),
  \end{equation}
  amongst all $\theta \in \Theta$.
\item[(c)] when a sequence of data sets $\{\Xi_T\}_{T=1}^\infty$ is
  generated by signals $(\tilde w, \tilde x_1, \ldots, \tilde x_N)$
  satisfying assumptions (A1),(A2) and (A4) for some $a_0 \in
  a_{\Theta_0}$, then for $\t_T = \sA(\Xi_T)$ the convergence of
  $a_{\theta_T}(x,w)$ to $a_0(x,w)$ takes place uniformly on every compact
  subset compact subset  of $\X \times \W$.
\end{itemize}

\section{Convex Parameterization of Models}
\label{sec:convex}
In this section we introduce the main construction of this paper: a
special class of stable projective parameterizations
$(\Theta,\Phi,\Psi)$, in which $\Theta$ are convex sets defined by a
family of linear matrix inequalities arrived at via an application of
the sum-of-squares relaxation, \cite{Parrilo03}.
The construction is motivated by
the earlier approaches from \cite{Bond10}, \cite{Tobenkin10}, and
\cite{Megretski08}, and is intended to improve consistency of the
associated system identification algorithms.

In the following definition, $x,\xi,\Delta$ and $q$ are 
real vector variables of dimensions $n_x$, and $w$ is a real vector
variable of dimension $n_w$.  In addition, $z=[\xi;x;w]$ and 
$v=[\xi;x;w;\Delta;q]$ are the real vector variables
of dimensions $2n_x+n_w$ and $4n_x+n_w$ respectively, 
constructed by concatenating $\xi$, $x$,
$w$, $\Delta$, and $q$.
Given an positive integer $N$ let
 \[  \mathcal P_N=\left\{p(z)=\sum_{\alpha \in \NN^{\dz} \; | \;
     \|\alpha\|_1\leq N}c_\alpha
 z^\alpha:~c_\alpha\in\mathbb R\right\}\]
 denote the set of all
 polynomials composed of monomials with scalar degrees no greater than $N$.

 Given a positive integer $N$, a positive constant $\d$, and
 a function 
$\Pi:~\mathbb R^{4n_x+n_w}\to\mathbb R^{n_\Pi}$,
let $\mathbf\Theta(N,\d,\Pi)$
be the set of all pairs $(\theta,r)$ of vectors
$\theta\in\mathbb R^{n_\theta}$ and $r \in {\cal P}_N$ for which there
exist matrices $P=P'\in\mathbb R^{n_x\times n_x}$,
$\Sigma_i=\Sigma_i'\in\mathbb R^{n_\Pi\times n_\Pi}$ (for
$i\in\{1,2\}$), and a positive scalar $\e$ such that
\begin{equation}\label{am3a}
P\geq \d I,~\Sigma_1\ge0,~\Sigma_2\ge0,
\end{equation}
\begin{flalign}\label{am3b}
r(z)+2\Delta'[e_\theta(x+\Delta)-e_\theta(x)]-|\Delta|^2_{P+C'C}&\\+|q|_P^2
-2q'(f_\theta(x+\Delta,w)-e_\theta(\xi))&=
\Pi(v)'\Sigma_1\Pi(v),\nonumber
\end{flalign}
\begin{flalign}\label{am3c}
2\Delta'[e_\theta(x+\Delta)-e_\theta(x)]-|\Delta|^2_{P+\e I}
&\\+|q|^2_P-2q'(f_\theta(x+\Delta,w)-f_\theta(x,w))&=
\Pi(v)'\Sigma_2\Pi(v)\nonumber,
\end{flalign}
where $e_\theta$ and $f_\theta$ are defined by \eqref{am2}.
By construction, $\mathbf\Theta(N,\d,\Pi)$ is a convex
set defined by a family of
linear matrix inequalities.
\begin{remark}
  The purpose of \eqref{am3b} is to
establish the condition
\begin{flalign} \label{eq:storage}
  r(z)+|e_\theta(x&+\Delta)- e_\theta(x)|^2_{P^{-1}} \\
  &\ge\nonumber\\
  |f_\theta(x+\Delta,w)&-e_\theta(\xi)|^2_{P^{-1}}+|C\Delta|^2,\nonumber
\end{flalign}
which in turn serves as a dissipation inequality used to bound
simulation error when using model \eqref{eq:dynamics} with $a=a_\theta$.
The purpose of \eqref{am3c} is to ensure that $e_\theta$ is a
bijection and 
establish the condition
\begin{flalign}\label{eq:lyap}
  |e_\theta(x+&\Delta)-e_\theta(x)|_{P^{-1}}^2\\
  &\ge\nonumber\\
  |f_\theta(x+\Delta,w)-&f_\theta(x,w)|^2_{P^{-1}}+\e|\Delta|^2, \nonumber
\end{flalign}
which is a non-linear version of the Lyapunov inequality, used
to prove that the model 
\eqref{eq:dynamics} with $a=a_\theta$ is $\ell^2$-incrementally stable.
\end{remark}


The following statement explains, partially, the utility of this
construction.

\begin{lem}\label{lem:delta}
If $\Theta$ is the set of all $\theta \in \RR^{n_\theta}$ such that
$(\theta,r) \in \mathbf\Theta(N,\d,\Pi)$ for some $r$,
then $(\Theta,\Phi,\Psi)$ is a stable projective parameterization.
Furthermore, for each $(\theta,r) \in
\mathbf\Theta(N,\d,\Pi) $ and data set $\Xi = (\tilde w,\tilde x_1, \ldots, \tilde x_n) \in \mathcal D(n_x,n_w,N,T)$ the function
\begin{equation}
\hat J_r(\Xi) :=  \frac{1}{NT} \sum_{t=1}^{N} \sum_{i=1}^{T} r(\tilde z_i(t))
\end{equation}
satisfies $\hat J_r(\Xi) \geq \frac{1}N \sum_{i=1}^N
J_T^{SE}(a_\theta,\tilde x_i(0),\tilde x_i, \tilde w)$.
\end{lem}

\begin{pf}
  For $\Theta$ as defined above to be a valid projective
  parameterization requires that $e_\theta$ be a bijection for all
  $\theta \in \Theta$.  The equality \eqref{am3c} holding for some
  $P \geq 0$ and $\Sigma_2 \geq 0$ implies that
  \begin{equation}
    2\Delta'(e_\theta(x+\Delta)-e_\theta(x)) \geq \e |\Delta|^2
  \end{equation}
 holds for all $x,\Delta \in \RR^n$.  As $e_\t$ is
  continuous, this condition implies $e_\t$ is a bijection
  (\cite{Poznyak08}, Theorem 18.15). 
  
  Next, we establish the connection between the conditions
  \eqref{am3b} and \eqref{am3c} and the inequalities 
  \eqref{eq:storage} and \eqref{eq:lyap}.
For all $a,b \in \RR^{n_x}$ and 
symmetric, positive definite $P \in \RR^{n_x\times n_x}$, the inequality
\begin{equation}
  |a|_{P^{-1}}^2 \geq 2b'a -|b|_P^2 \label{eq:restrict}
\end{equation}
holds due to the fact that $|a-Pb|_{P^{-1}}^2 \geq 0$.
Fixing any $(\theta,r) \in \mathbf\Theta(N,\d,\Pi)$, and
applying \eqref{eq:restrict} with $a = e_\theta(x+\Delta)-e_\theta(x)$ and $b = \Delta$, we see that
there exists a $P \in \RR^{\dx\times\dx}$ such that:
\begin{flalign*}
  r(z)+|e_\theta(x+\Delta)-e_\theta(x)|^2_{P^{-1}}-|C\Delta|^2&\\
  +|q|_P^2
-2q'(f_\theta(x+\Delta,w)-e_\theta(\xi)) &\geq 0,
\end{flalign*}
and
\begin{flalign*}
|e_\theta(x+\Delta)-e_\theta(x)|^2_{P^{-1}}-\e|\Delta|^2&\\
+|q|^2_P-2q'(f_\theta(x+\Delta,w)-f_\theta(x,w)) &\geq 0,
\end{flalign*}
hold for all $x,\xi,\Delta,q$ in $\RR^{n_x}$ and $w \in \RR^{n_w}$.
Analytically minimizing these expressions with respect to $q$
demonstrates these inequalities imply \eqref{eq:storage} and
\eqref{eq:lyap} hold for all $x,\xi,\Delta \in \RR^{n_x}$ and $w \in
\RR^{n_w}$.

Fix $x_{01}, x_{02} \in \RR^{n_x}$ and $w: \ZZ_+ \mapsto \RR^{n_w}$,
and let $x_i$ be the solution $x_i = G_{a_\theta}(x_{0i},w)$ for $i \in \{1,2\}$. The inequality
\eqref{eq:lyap} with $x \equiv x_1(t)$ and $\Delta \equiv
x_2(t)-x_1(t)$ implies:
\begin{flalign*} |e_\theta(x_{02})-e_\theta(x_{01})|_{P^{-1}}^2 \geq&\:
  |e_\theta(x_2(T))-e_\theta(x_{1}(T))|_{P^{-1}}^2\\
  &\:+\e\sum_{t=0}^{T-1}
|x_1(t)-x_2(t)|^2.
\end{flalign*}
As $P \geq 0$, we conclude that \eqref{eq:dynamics} with $a\equiv
a_\theta$ is $\ell^2$-incrementally stable.
Take $x_i = G_{a_\theta}(\tilde
x_i(0),\tilde w)$, then examining \eqref{eq:storage} with $z = \tilde
z_i(t)$ and $\Delta = x_i(t) - \tilde x_i(t)$ leads to
$$ \sum_{t=1}^{T} r(\tilde z_i(t)) \geq
|e_\theta(x_i(T))-e_\theta(\tilde x_{i}(T))|_{P^{-1}}^2+\sum_{t=0}^{T-1}
|C(x_i(t)-\tilde x_i(t))|^2,$$
from which one can readily conclude $\hat J_r(\Xi) \geq J^{SE}_T(a_\theta,\Xi)$.
\end{pf}

The following definition provides a subset of systems for which we can establish consistency results.

\begin{dfn}
  The {\it recoverable set} defined by $\mathbf
  \Theta(N,\d,\Pi)$ is the set of functions $a:
  \RR^{\dx}\times\RR^{\dw} \mapsto \RR^{\dx}$ such that there exists a
  pair $(\theta,r) \in \mathbf \Theta(N,\d,\Pi)$ with
  $a(x,w) \equiv a_\theta(x,w)$ and
  \begin{equation}
  r([a(x,w);x;w]) = 0 \label{eq:requality}
  \end{equation}
  for all $(x,w) \in \RR^{\dx}\times\RR^{\dw}$.
\end{dfn}

The following lemma establishes that the subset of $a_\Theta$
consisting of recoverable models can be made to include all stable, linear
state-space models of appropriate dimensionality. 

\begin{lem}\label{lem:recovery}
 Let $\Phi$ and $\Psi$ be
  finite sequences of real analytic functions, as above, whose
  respective spans include all linear functions, and let $\Pi:
  \RR^{4\dx+\dw}\mapsto \RR^{n_\Pi}$ be a function such that the span of its components
  include all linear functions.  Then for all $C \in
  \RR^{\dy\times\dx}$ and $N \geq 2$ the recoverable set defined by $\mathbf
  \Theta(N,\d,\Pi)$ includes all functions $a: \RR^{\dx}
  \times \RR^{\dw} \mapsto \RR^{\dx}$ given by
   \begin{flalign}
    a(x,w) = Ax + Bw,
  \end{flalign}
  where $A \in \RR^{\dx \times \dx}$ is Schur (stable) and $B \in
  \RR^{\dx \times\dw}$.
\end{lem}

\begin{pf}
  As $A$ is Schur, there exists a symmetric positive definite matrix
  $P$ solving the Lyapunov equation $P-A'PA=C'C+\d I$ for any positive
  $\d$.  Choose $\t$
  such that:
  $$e_\t(x) = Px, \quad f_\t(x,w) = Pa(x,w).$$
  The constraint \eqref{am3b} is therefore equivalent to
  \begin{flalign*}
    r(z) + \D'P\D +|q|^2_{P}  &\\
   +2q'(PA\D-e(\xi)+f(x,w))-|C\D|^2 &=
    \Pi(v)'\Sigma_1\Pi(v).
  \end{flalign*}
Explicit minimization w.r.t. $q$  shows that the polynomial on the
left hand side of this equality is lower bounded by:
 $$r(z) + \D'P\D -|PA\D-e(\xi)+f(x,w)|_{P^{-1}}^2-|C\D|^2.$$
  This is a concave function of $\D$ as $\D'(P-A'PA-C'C)\D
  =-\d|\D|^2$.  Explicit minimization w.r.t $\D$ provides
  $$r(z) -|e(\xi)-f(x,w)|_Q^2$$
  as a lower bound for the original polynomial for some $Q = Q' \geq
  0$.  The function
  $r(z) = |e(\xi)-f(x,w)|_Q^2$,
  belongs to ${\cal P}_N$ and clearly satisfies
  \eqref{eq:requality}.
  Furthermore, with this choice of $r(z)$ the left hand side of \eqref{am3b} is a non-negative
  quadratic polynomial so that there exists
  an appropriate choice of $\Sigma_1 \geq 0$ to ensure \eqref{am3b} holds.
  A similar analysis shows that an appropriate choice of $\Sigma_2
  \geq 0$ also exists, thus $a$ belongs to the recoverable set.
\end{pf}

A simple example of a nonlinear function  $a: \RR \times \RR \to \RR$ belonging to such a recoverable
set is given by
\[ e(a(x,w)) = \frac{1}{2} x + b(w) \]
where $e(x) = \frac{3}2 x + x^3$ and $b: \RR \to \RR$ is an arbitrary
polynomial.   That an appropriate
recoverable set exists is shown in the appendix.

\section{Identification Algorithm}
\label{sec:noise}
This section presents an algorithm for transforming data sets 
$\Xi = (\tilde w, \tilde x_1, \ldots, \tilde x_N) \in {\cal D}(\dx,\dw,N)$ into parameter
vectors $\hat \t \in \RR^{n_\theta}$, followed by an asymptotic analysis of the algorithm.
For the remainder of this section we define
\[
\aleph = \{ \alpha \in \NN^{\dz} \; | \; 2\|\alpha\|_1 \leq N\}.
\]
\subsubsection*{Algorithm $\sA(\d,\Pi,\kappa)$}

\begin{enumerate}
\item[(i)] Select a constant $\d >0$ and a function $\Pi: \RR^{4\dx+\dw} \to \RR^{n_\Pi}$,
  as described in Section~\ref{sec:convex}. Additionally, select a 
  constant $\kappa \in (0,\infty]$.
\item[(ii)]
  Form the matrix $\tilde M_\Xi \in
  \RR^{|\aleph|\times|\aleph|}$ given by:
  $$
  [\tilde M_\Xi]_{j_1,j_2} = \tilde\mu_{\alpha_{j_1}+\alpha_{j_2}}(\Xi),
  $$
  where $\tilde \mu_\alpha(\cdot)$ are the linearized empirical moments defined by \eqref{eq:tildemu}, and let $\hat
  M_\Xi$ be the projection of $\frac{1}2 (\tilde M_\Xi + \tilde M_\Xi')$ onto the
  closed convex cone of positive semidefinite matrices.
\item[(iii)]
  Find the $\t \in \RR^{n_\t}$, $r \in {\cal P}_N$ and $R = R'
  \in \RR^{|\aleph|\times|\aleph|}$ that minimize:
  \[ \trace(R\hat M_\Xi) \]
  subject to $(\t,r) \in \mathbf
  \Theta(N,\d,\Pi),$
  \[ r(z) = \sum_{j_1=1}^{|\aleph|}
  \sum_{j_2=1}^{|\aleph|} [R]_{j_1,j_2}
  z^{\alpha_{j_1}+\alpha_{j_2}},\]
  and $\|R\|_F^2 \leq \kappa$.  Take $(\hat \theta,\hat r,\hat R)$ to
  be an optimizing $(\theta, r, R)$.
\end{enumerate}
\begin{remark}
  Note that the algorithm is well defined if any subset of the set of
  vector degrees $\aleph$ is substituted in lieu of $\aleph$.
\end{remark}

\begin{remark}
 Examining the definition of $R$ and $\hat M$, one sees that when
 $\tilde x_1 = \tilde x_2 = \ldots = \tilde x_N$, $\tilde M_\Xi = \hat
 M_\Xi$ and the objective function
 $\trace(R\hat M_\Xi)$ is equal to $\hat J_{\hat r}(\Xi)$, the previously
 established upper bound on simulation error.  
The additional parameter $\kappa$, when finite, ensures that the
optimal value of the optimization problem of step (iii) has a
continuous dependence on $\hat M_\Xi$ (and by extension, the linearized
empirical moments).  
\end{remark}
\subsection{Asymptotic Analysis}
This section analyzes the properties of algorithm $\sA$ when data
sets are generated according to one of the two data generation
mechanisms described in Section~\ref{sec:mechanism}.  By $(\tilde w,
\tilde x_1, \ldots, \tilde x_N)$, $(v_1, \ldots, v_N)$, and $(\bar
x_1, \ldots, \bar x_N)$, we mean those signals described in
assumptions (A3) or (A4).  Let $\tilde w^{(T)}$, $\tilde
x^{(T)}$ and $\bar x_i^{(T)}$ are the restrictions of $\tilde w,
\tilde x_i$ and $\bar x_i$ to $\{0,\ldots,T\}$. We define $
\bar M_T \in \RR^{|\aleph|\times|\aleph|}$ to be a matrix of ``noiseless''
empirical moments, given by
$$ [\bar M_T]_{j_1,j_2} = \frac{1}{N} \sum_{i=1}^N
\mu_{\alpha_{j_1}+\alpha_{j_2}}(\bar x_i^{(T)},\tilde w^{(T)}).$$
The following lemma demonstrates that the
linearized empirical moments, under suitable assumptions, converge to
these noiseless empirical moments.

\begin{lem}\label{lem:momapprox}
 Let $\tilde w: \NN \to \W$ and $\tilde x_i: \NN \to \X$, for $i \in
  \{1,\ldots,N\}$, satisfy  assumptions (A1)-(A3).  Then
  $$\hat M_{\Xi_T} - \bar M_T \cvgip 0$$
  as $T \rightarrow \infty$.
\end{lem}

The following statement justifies the use of this algorithm
under the assumptions (A1)-(A3). 
\begin{thm}\label{thm:simerr}
  Let $\tilde w: \NN \to \W$ and $\tilde x_i: \NN \to \X$, for $i \in
  \{1,\ldots,N\}$, satisfy assumptions (A1)-(A3). 
  Fix $\d$, $\Pi$, and $\kappa$ as in algorithm
  $\sA$, and  let $(\t_T,r_T,R_T)$ be the $(\hat \t, \hat r,\hat R)$ found by
  applying algorithm $\sA(\d,\Pi,\kappa)$ to the data set $\Xi_T$.  If
  $\kappa$ is finite, then for every $\e > 0$ and $\gamma \in
  (0,1)$ there exists a $T$ such that, with probability $1-\gamma$,
  $$
 \e + \trace(R\bar M_T) \geq \trace(R_T\bar M_T ) 
  \geq \frac{1}{N} \sum_{i=1}^N \SE_T(a_{\t_T},\bar x_i(0),\bar
  x_i^{(T)},\tilde w^{(T)}),
  $$
  holds for all  $(\t,r,R)$ feasible for the optimization problem in step
  (iii) of $\sA(\d,\Pi,\kappa)$.
\end{thm}
The proof of the above theorem is in the appendix.

The following result in characterizes the  consistency of
$\sA$ in terms of recoverable functions.

\begin{thm}\label{thm:consistency}
  Fix $\Pi: \RR^{4\dx+\dw}\to\RR^{n_\Pi}$, $\kappa \in (0,\infty)$,
  $\d > 0$.  Let $\{\Xi_T\}_{T=1}^\infty$ be a sequence of data sets
  defined by signals $\tilde w: \NN \to \W$ and
  $\tilde x_i: \NN \to \X$, for $i \in \{1,\ldots,N\}$, 
  satisfying assumptions (A1), (A2), and (A4), where the function
  $a_0: \X \times \W \to \X$ described in assumption (A4) is in the recoverable set
  defined by $\mathbf \Theta(N,\d,\Pi)$.
 If $\t_T$ is the parameter vector found by applying algorithm $\sA(\d,\Pi,\kappa)$ to
  $\Xi_T$, then
$$\sup_{(x,w) \in K} \{ |e_{\t_T}^{-1}(f_{\t_T}(x,w)) -
a_0(x,w)|\} \cvgip 0$$
as $T \rightarrow \infty$ for all  compact sets $K \subset \X
  \times \W$.
\end{thm}
The proof of this theorem and supporting lemmas are given in the
appendix.

\section{Example}
\label{sec:example}
This section examines the performance of algorithm $\sA$ via a simple simulated example. Two alternative identification algorithms are introduced based on least square minimization; this is followed by a comparison of these algorithms and $\sA$ on a simulated data set.
\subsection{Least Squares Identification Approaches}
Let $\Phi = \{\phi_i\}_{i=1}^{n_\t}$ and $\Psi = \{\psi_i\}_{i=1}^{n_\t}$ be fixed sequences of polynomial functions, $\phi_i: \X \times \W \to \X$ and $\psi_i: \X \to \X$.  Let $x$ and $\D$ denote real vector variables of dimension $\dx$.
Given a positive constant $\d$ and function $\Lambda: \RR^{2\dx} \to \RR^{n_{\Lambda}}$, let $\mathbf \Omega(\d,\Lambda)$ be the set of $\t \in \RR^{n_\t}$ such that there exists some $\Sigma \in \RR^{n_{\Lambda}\times n_{\Lambda}}$ such that:
\[
  \Sigma \geq 0,
\]
\[
 2\D'(e_\t(x+\D)-e_\t(x)) - \d |\D|^2 = \Lambda([x;\D])'\Sigma \Lambda([x;\D]).
\]
From the proof of Lemma~\ref{lem:delta} we see that $\t \in \Omega(\d,\Lambda)$ guarantees $e_\t$ is a bijection.

The following two algorithms produce a parameter vector $\hat \t$ from a data set $\Xi = (\tilde w, \tilde x_1, \ldots, \tilde x_N) \in {\cal D}(\dx,\dw,N,T)$.
\subsubsection{Least Squares Algorithm}
Take $\hat \t$ to be a $\t \in \Omega(\d,\Lambda)$ that minimizes
$$ \frac{1}{NT} \sum_{i=1}^N \sum_{t=1}^T |e_\t(\tilde x_i(t)) - f_\t(\tilde x_i(t-1),\tilde w(t-1))|^2.$$

The following algorithm adapts the least squares objective to use the linearized empirical moments $\tilde \mu_\alpha(\cdot)$ defined by \eqref{eq:tildemu} with the aim of bias elimination.
\subsubsection{Modified Least Squares Algorithm}
\begin{enumerate}
\item[(i)]
  Fix a $\d > 0$ and $\kappa \in (0,\infty]$, and let $\aleph = \{\alpha_j\}_{j=1}^{|\aleph}\subset \NN^{\dz}$ be the smallest set of vector degrees such that for each $\t \in \Theta$ there exists coefficients $c_\alpha \in \RR$ satisfying $e_\t(\xi) - f_\t(x,w) = \sum_{\alpha\in\aleph } c_\alpha [\xi;x;w]^\alpha.$
\item[(ii)] Define $\tilde M \in \RR^{|\aleph|\times|\aleph|}$ by $[\tilde M]_{j_1,j_2} = \tilde \mu_{\alpha_{j_1}+\alpha_{j_2}}(\Xi)$, and take $\hat M$ to be the projection of $\frac{1}2(\tilde M + \tilde M')$ on the cone of positive semidefinite matrices.
\item[(iii)]  Find the $\t \in \Omega(\d,\Lambda)$, satisfying $|\theta|^2 \leq \kappa$, that minimizes
  $$\trace\left(\t\t' \sum_{i=1}^{n_x} \Gamma_i'\hat M\Gamma_i\right),$$
  where $\Gamma_i: \RR^{|\aleph|\times n_\t}$ is defined so that $[e_\t(\xi)-f_\t(x,w)]_i = \sum_{\alpha \in \aleph} z^{\alpha_j} [\Gamma_i \t]_j$.
\end{enumerate}

\subsection{Simulated Example}

We provide a simple simulated example to compare the performance algorithm to the least squares algorithms defined above.
We examine a SISO nonlinear output error data set generated in the following fashion.  The input is a scalar sequence $\tilde u: \{0,\ldots,800\} \mapsto \RR$, generated as i.i.d. random variables distributed uniformly on $[-1,1]$.
We examine the response of the system:
$$ \bar x(t+1) + \frac{1}5 \bar x(t+1)^5 = \frac{1}3 \bar x(t)^3+ 5\tilde u(t)$$
starting from $\bar x(0) = 0$.  For $i \in \{1,\ldots,10\},$ observed data sequences, $\tilde x_i(t)$, are generated according to $\tilde x_i(t) = \bar x(t) + v_i(t)$ where the $v_i(t)$ are i.i.d. zero mean Gaussian random variables with variance 0.09 and independent across trials and from the input, leading to a  signal-to-noise ratio of approximately $25$ dB.  We take $C = 1$, $\Phi$ to contain all monomials of degree less than or equal to five and $\Psi$ to contain all monomials affine in $u$ and of degree less than or equal to three in $x$.

The identified models are computed on a subset of the available data  revealing only the samples with $t \in \{0,\ldots,T_h\}$ for each  $T_h = 100\cdot 2^{h}$ for $h \in \{0,1,2,3\}$ .
The parameters $\d$ and $\kappa$ were taken to be $0.01$ and $\infty$ respectively, and the set $\aleph$ from the definition of the modified least squares algorithm was also used for algorithm $\sA$.  The sum-of-squares programs were prepared using YALMIP \cite{Lofberg09}.  The choices of $\Pi$, as in algorithm $\sA$, and $\Lambda$, as in the modified least squares algorithm, that this software makes ensure that $\Theta(N,\d,\Pi) \subset \Omega(\d,\Lambda)$, i.e. the modified least squares algorithm searches over a larger set of models than algorithm $\sA$.

\begin{figure*}
  \centering
  \includegraphics[width=\textwidth]{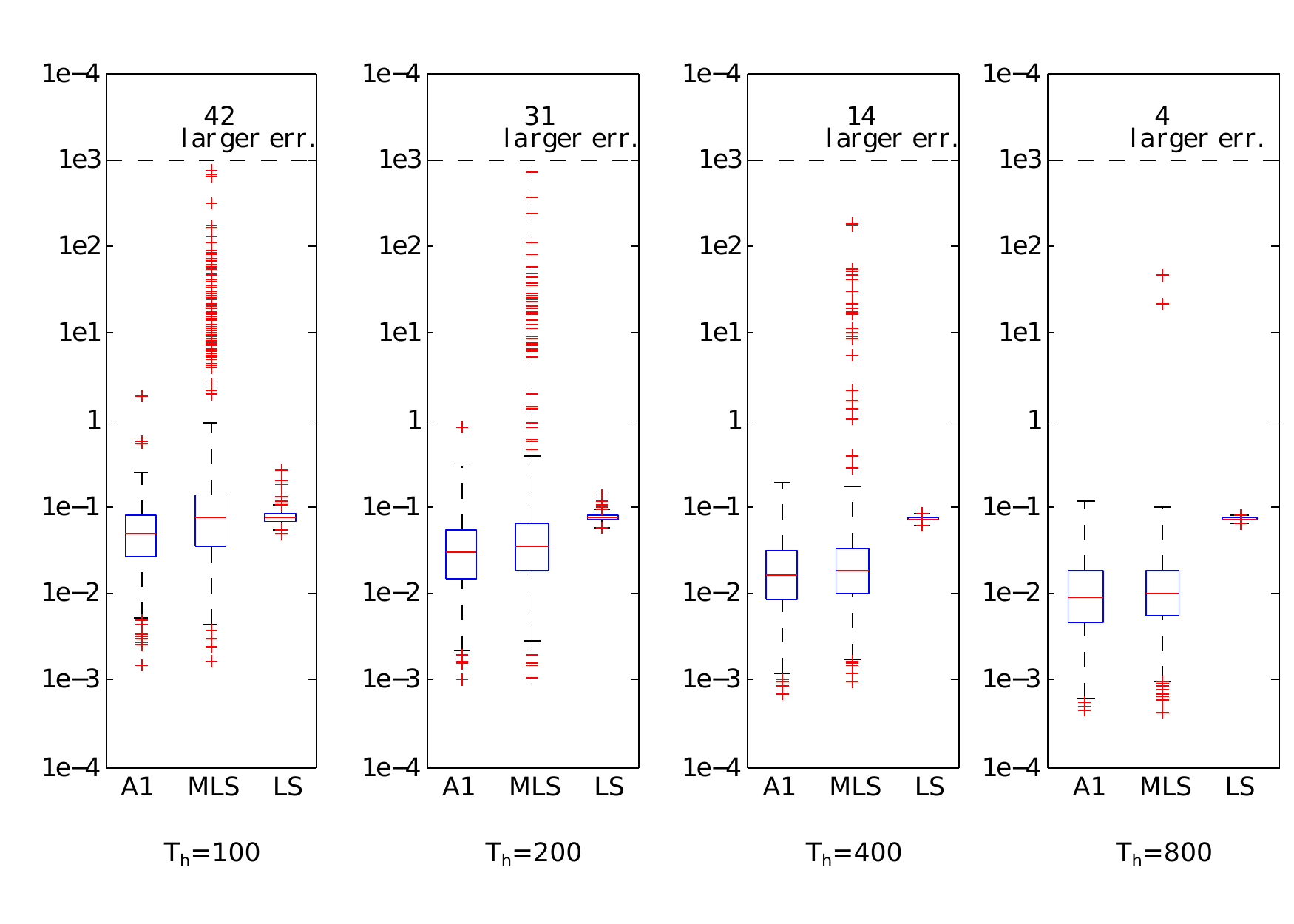}
  \caption{
    Comparison of the algorithm $\sA$ (A), the modified least squares algorithm (MLS) and the least squares algorithm (LS) for various training horizons.  Plotted on a log scale is the distribution of the normalized simulation error for the validation input over 1,000 realizations. $T_h$ indicates the number of training samples.
    The vertical scale of $1$ indicates identical performance to the model that always outputs zero (i.e. $100$ percent simulation error) and $1e-2$ indicates 1 percent simulation error.
    At the top of each plot is the number of MLS models having greater than $1000$ percent simulation error.
  }
  \label{fig:nl_iv_1st_d3}
\end{figure*}

To validate the models we generate an additional input sequence $\bar u: \{0,\ldots,800\} \mapsto \RR$, again i.i.d. uniformly distributed on [-1,1], that is independent from the training input, and noise.
We compute the response $\bar x_{\textrm{val}}(t)$ of the true system to this input from zero initial conditions and compute a normalized simulation error:
$$  \sum_{t=0}^{T_h} | \bar x_{\textrm{val}}(t) - x_{h}(t)|^2 \Big / \sum_{t=0}^{T_h} | \bar x_{\textrm{val}}(t)|^2$$
where $x_{h}(t)$ is the response of the optimized model to the same input and starting from the same initial condition.  These calculations were performed for 1,000 independent realizations of the problem.  Figure~\ref{fig:nl_iv_1st_d3} plots a comparison of the models generated by the three algorithms

As the amount of available data increases, the distribution of validation simulation errors tends toward zero for both algorithm $\sA$ and the modified least squares approach. By contrast the result of the least squares approach without modification remains biased, though the variance of errors decreases.   One sees that the modified least squares algorithm  generates a large number of poorly performing models and generally under-performs algorithm $\sA$ in terms of median as well.  Note that the vertical scale in these plots is logarithmic: at $T_h = 200$ the majority of the  models rendered by algorithm $\sA$ have less that $3$ percent validation simulation error.


\section{Conclusions}
\label{sec:conc}

In this paper we have presented a convex optimization approach to
nonlinear system identification from noise-corrupted data. The main
contributions are a particular convex family of stable nonlinear
models, and a procedure for identifying models from data based on
empirical moments.

This builds upon previous work by the authors \cite{Tobenkin10}, \cite{Bond10} and \cite{Megretski08} and offers two main advantages
over the previously proposed methods: the complexity of the computation does not
grow with the length of the data, and the empirical moments can be
``mixed'' from different experiments to achieve a consistent estimator
when the true system is in the model class. This is reminiscent of the
instrumental variables approach to total least squares, although we
suggest minimizing an alternative convex criterion. The advantages of
the proposed method over least squares methods were illustrated via a
simple simulated example.

\appendix[Proofs]
\subsection{Modeling Scenario Satisfying (A1) and (A4)}
The section provides an example of a class of dynamical systems and inputs that
result in signals satisfying conditions (A1) and (A4).
\begin{prop}
  Let $\tilde w(t) \in \W$ be a vector valued i.i.d. stochastic
  process with $\tilde w(t)$ having an absolutely continuous
  distribution with a lower semicontinuous density whose support is
  bounded and contains $0$.  Let $a_0: \X \times \W$ be a real
  analytic function such that \eqref{eq:dynamics} with $a = a_0$ is
  BIBO and incrementally exponentially stable and the linearization of
  $a_0$ about the unique equilibrium $(x,u) = (x_0,0)$, with $x_0 =
  a_0(x_0,0)$, is controllable.  Then for any $\{\bar x_{i0}\}_{i=1}^N
  \subset \X$ independent of $\tilde w(t)$, the signal $\tilde w(t)$
  and the signals $\bar x_i: \NN \to \X$, each defined by
  \eqref{eq:dynamics} with $w(t) \equiv \tilde w(t)$ and  $\bar x_i(t)
  \equiv x(t)$ with $x(-1) = \bar x_{i0}$, satisfy (A1) and (A4).
\end{prop}
\begin{pf}
  That $\tilde w(t)$ satisfies (A1) is immediate.
 We use several results from \cite{Meyn91}.  The controllability
  assumptions given imply that the Markovian system generating the
  sequences $(\bar x_i(t),\tilde w(t))$ is
  weakly stochastically controllable (w.s.c) (\cite{Meyn91}, Corollary
  2.2).  The BIBO stability of \eqref{eq:dynamics} with $a=a_0$
  and boundedness of $\tilde w(t)$ ensure $(\bar x_i(t), \tilde
  w(t))$ is bounded in probability. This
  boundedness in probability, the w.s.c. of the system, and the conditions on the
  distribution of $\tilde w(t)$ imply that the Markovian
  system generating $(\bar x_i(t), \tilde w(t))$ is a
  positive Harris recurrent Markov Chain.  Thus there exists a
  unique stationary probability distribution $\pi$ that is independent
  of initial condition.  The w.s.c of the system and the assumptions
  on the distribution of $\tilde w(t)$ ensure that the support of
  $\pi$ contains an open set.  As the support of $\pi$ is clearly bounded
  \begin{equation}
    \lim_{T\rightarrow\infty}\frac{1}{T} \sum_{t=0}^{T-1} h(\bar
    x_i(t),\tilde w(t)) = \int h d\pi,
  \end{equation}
  for all continuous $h: \X \times \W \to \X$ (\cite{Meyn91} Proposition 2.3).
  From this we conclude that each pair $(\tilde w, \bar x_i)$ is
  persistently excited with respect to $a_0$, as in Definition~\ref{def:persist}.
\end{pf}

\subsection{A Simple Recoverable Nonlinear System}
Let $C = 1$ and define $a: \RR \times \RR \mapsto \RR$ by:
\[ e(a(x,w)) = \frac{1}{2} x + b(w) \]
where $e(x) = \frac{3}2 x + x^3$ and $b: \RR \to \RR$ is a
polynomial.  Let $\Phi$ and $\Psi$ be sequences of real analytic
functions as in the definition of $e_\theta$ and $f_\theta$ such that
there exists a $\theta$ with $e = e_\theta $ and $f = f_\theta$.
Let $x,\xi,\D,q,w$ be real numbers, $z = [\xi;x;w]$ and $v = [\xi;x;w;\D;q]$.
With $r(z) = 2(e(\xi)-f(x,w))^2$ and $P = 1$, clearly
\eqref{eq:requality} is satisfied and the left hand side of the
equality \eqref{am3b} can be expressed as
\begin{flalign*} (\sqrt{2}(e(\xi)-f(x,w)) -q/\sqrt{2})^2 + (\D -q/2)^2 \\+
  \frac{1}4q^2+ (\sqrt{6} x\D + \sqrt{3/2}\D^2)^2 + \frac{1}2\D^4,
\end{flalign*}
i.e. as the sum of squares of polynomials.  Similarly the left hand
size of \eqref{am3c} can be expressed as:
\[ (\D + q/2)^2 + \frac{3}4q^2+
(\sqrt{6} x\D + \sqrt{3/2}\D^2)^2 + \frac{1}2\D^4.\]
These sum of squares decompositions show that there is a choice of
$\Pi(v)$ such that matrices $\Sigma_1 \geq 0$ and $\Sigma_2 \geq 0$ 
satisfying \eqref{am3b} and \eqref{am3c} are guaranteed to exist.
For this $\Pi$ and $N \geq 6$, the polynomial $r$ belongs to $ {\cal
  P}_N$ so that $a$ is in the recoverable set defined by $\mathbf \Theta(N,1,\Pi)$.

\subsection{Proof of Lemma~\ref{lem:momapprox}}
\begin{pf}
  By assumptions (A1) and (A3),  the matrices $\tilde M_{\Xi_T}$
  are uniformily bounded in $T$ with probability one. 
  Since each $\bar M_T \geq 0$ and projection onto the positive semidefinite cone is a continuous
  function, this boundedness implies it is sufficient to show that 
 $$\tilde \mu_\alpha(\Xi_T) - \mu_{\alpha}(\bar x_i^{(T)},\tilde
  w^{(T)})\cvgip 0,$$
  as $T \rightarrow \infty$, for each $i \in \{1,\ldots, N\}$ and
  $\alpha \in \NN^{\dz}$ with $\|\alpha\|_1 \leq N$.
 Let $\bar z_i(t) = [\bar x_i(t);\bar x_i(t-1);\tilde w(t-1)]$ and $\zeta_i(t) = \tilde z_i(t) -
  \bar z_i(t)$. For all $\alpha \in \NN^{\dz}$ with $\|\alpha\|_1 \leq
  N$ we have:
  \begin{equation}
    \prod_{i=1}^{\|\alpha\|_1} [\tilde z_i(t)]_{\beta_{\alpha}(i)}
    =\sum_\kappa \nu_\kappa(t),
  \end{equation}
  where the sum runs over all $\kappa: \{1,\ldots,\|\alpha\|_1\} \to
  \{0,1\}$ and $\nu_\kappa(t)$ is defined by:
  \begin{equation}
    \nu_\kappa(t) := \prod_{i=1}^{\|\alpha\|_1} [\bar z_i(t)]_{\beta_{\alpha}(i)}^{(1-\kappa(i))}[\zeta_i(t)]_{\beta_{\alpha}(i)}^{\kappa(i)}.
  \end{equation}
  We see:
 \begin{equation}
   \frac{1}{T}\sum_{t=1}^{T} \nu_{0}(t) = \frac{1}{T}\sum_{t=1}^{T} \prod_{i=1}^{\|\alpha\|_1} [\bar
    z_i(t)]_{\beta_{\alpha}(i)}.
  \end{equation}
  Condition \eqref{eq:yconv} of (A3) implies that $\lim_{t\rightarrow \infty}|\bar z_i(t)-\bar z_j(t)| = 0$ which, combined with the boundedness of $\bar z_i(t)$, allows us to conclude:
  \begin{equation}
    \lim_{T\rightarrow\infty} \left |\frac{1}{T}\sum_{t=1}^{T} \nu_{0}(t)-\mu_{\alpha}(\bar
      x_i^{(T)},\tilde w^{(T)})\right|  = 0.
  \end{equation}

  Each $\nu_\kappa(t)$, for $\kappa \neq 0$, is obtained as a multi-linear
  function of the noise sequences $v_i(t)$, whose coefficients depend
  only on $\{\bar z_i(t)\}_{i=1}^N$.  As the $v_i(t)$ are zero mean, bounded, and
  independent from one another and each $\bar z_i(t)$, we see $\nu_\kappa(t)$ is a zero mean, bounded random
  variable. In addition, condition \eqref{eq:yconv} of (A3) and the boundedness of $\bar x_i(t)$ imply the correlation between
  $\nu_\kappa(t)$ and $\nu_\kappa(\tau)$ decays exponentially
  as $|t-\tau|\rightarrow \infty $, hence $\frac{1}{T}
  \sum_{t=1}^{T} \nu_\kappa(t) \cvgip 0$.
\end{pf}

\subsection{Proof of Theorem~\ref{thm:simerr}}
\begin{pf}
  Optimality of $(\t_T,r_T,R_T)$ implies that:
  \begin{equation}
     \trace(R(\hat M_{\Xi_T}-\bar M_T))+\trace(R_T (\bar M_T-\hat
     M_{\Xi_T}))+\trace( R\bar M_T)\geq \trace(R_T \bar M_T)
   \end{equation}
   for all feasible $(\t,r,R)$. From this, and the fact that $\kappa
   \geq \|R\|_F^2$ for all feasible $R$, we see:
 \begin{equation}
     2\sqrt{\kappa}\|\bar M_T - \hat M_T\|_F+\trace( \bar M_T R)\geq
     \trace(\bar M_T R_T).
   \end{equation}
   Fixing $\gamma \in (0,1)$, Lemma~\ref{lem:momapprox} guarantees
   there exists a $T$ sufficiently large such that $\|\bar M_T-\hat M_T\|_F \leq
   \frac{\e}{2\sqrt{\kappa}}$ with probability $1-\gamma$.  Taking this
     $T$, the result follows by examining the definition of
     $r(\cdot)$, $\bar M_T$ and Lemma~\ref{lem:delta}.
\end{pf}

\subsection{Proof of Theorem~\ref{thm:consistency} and Supporting Lemmas}
We begin by proving several supporting Lemmas.

\begin{lem}\label{lem:a}
  For any feasible point $(\t,r,R)$ of the optimization of step (iii)
  in algorithm $\sA$, the function $a_\t =
  e_\t^{-1}\circ f_\t$ is real analytic.  Furthermore, we have:
  \begin{equation}
    r(z) \geq \frac{1}{\d}|e_\t(\xi)-f_\t(x,w)|^2 \geq \frac{\d}4|\xi - a_\t(x,w)|^2 ,\label{eq:abound}
  \end{equation}
  where $z = [\xi ; x ; w]$, with
  $\xi,x \in \X$ and $w \in \W$.
\end{lem}
\begin{pf}
  Examining \eqref{am3c} with $x = x_1$ and $\D = x_1-x_2$
  we see that:
  \begin{equation}
    2\D'(e_\t(x+\D)-e_\t(x)) \geq |\D|^2_P 
  \end{equation}
  for all $x,\D \in \X$.  Letting $E(x) = \frac{\partial e_\t}{\partial
    x}$, the above implies $E(x) + E(x)' \geq  P \geq \d I$ for
  all $x \in \X$ as $e_\t(\cdot)$ is continuously differentiable.  Thus
  $\det(E(x)) \neq 0$ for all $x \in \X$, and as $e_\t(\xi)-f_\t(x,w)$ is
  real analytic the Inverse Function Theorem for
  holomorphic maps
  (\cite{Fritzsche02} Chapter I, Theorem 7.6) ensures $e_\t^{-1}(f_\t(x,w))$ 
  is real analytic.

  The condition \eqref{am3b} for  $x_1 = x_2 = x$ gives
 \begin{equation}
    r(z) \geq |e_\t(\xi)-f_\t(x,w)|^2_{P^{-1}} \geq \frac{1}{\d}|e_\t(\xi)-f_\t(x,w)|^2,
  \end{equation}
  as $P \geq \d I$.
Furthermore, 
  \begin{flalign*}
    | e_\t(\xi) - f_\t(x,w)|^2 =&\: |e_\t(\xi) - e_\t(a_\t(x,w))|^2\\
    =&\: \left |\xi-a_\t(x,w)\right|^2_{\bar E'\bar E},
  \end{flalign*}
  where $\bar E = \int_0^1 E(a_\t(x,w)+\tau(\xi-a_\t(x,w)))d\tau$.  We see
  that
  $$\bar E'\bar E \geq \frac{\d}2 (\bar E + \bar E') -\frac{\d^2}4 I \geq \frac{\d^2}4 I,$$
  where the first inequality follows from $(\bar E-\frac{\d}2 I)'(\bar E-\frac{\d}2 I)
  \geq 0$
  and the second inequality from $\bar E$ being a convex combination
  of point evaluations of $E(x)$.  From this we can conclude \eqref{eq:abound}
\end{pf}

Now we present the proof of Theorem~\ref{thm:consistency}.
\begin{pf}
  Fix any compact set $K \subset \X \times \W$.  By assumption there
  exists a $(\t_0,r_0,R_0)$, feasible for the optimization in step
  (iii) of algorithm $\sA(\delta,\Pi,\kappa)$ such that $a_0 =
  a_{\t_0}$ and $ r_0(a_0(x,w),x,w) \equiv 0$. This implies
  $\trace(R_0\bar M_T) = 0$, so that, by Lemma~\ref{lem:momapprox}, 
  $\trace(R_0\hat M_{\Xi_T}) \cvgip 0$.

  As $(\t_0,r_0,R_0)$ is in the feasible set of the optimization in
  step (iii) of $\sA(\delta,\Pi,\kappa)$ and 
 $(\t_T,r_T,R_T)$ is optimal,
  \begin{equation}
   \trace(R_0\hat M_{\Xi_T}) \geq \trace(R_T \hat M_{\Xi_T}),
 \end{equation}
 so that $\trace(R_T\hat M_{\Xi_T}) \cvgip 0$ as well.  Moreover,  as
 $\|R_T\|_F^2 \leq \kappa$,
 \begin{equation}
  \sqrt{\kappa}\|\hat M_{\Xi_T}-\bar M_T\|_F +\trace(R_T \hat
  M_{\Xi_T}) \geq \trace(R_T  \bar M_T),
\end{equation}
by Cauchy-Scharwz.
Lemma~\ref{lem:momapprox} now implies $\trace(R_T  \bar M_T)  \cvgip 0$.

As each $r_T$ is a polynomial and the matrices $R_T$ are uniformily
bounded,

As each $(\bar x_i,\tilde w)$ is
persistently exciting with respect to $a_0$
(Definition~\ref{def:persist}),  there exists of a positive
measure $\pi$ on the space $\X\times \W$, supported on an open set,
such that for all $\e > 0$ there exists a $T_0 \in \NN$ with
$$ \e\|R_T\|_F + \trace(R_T \bar M_T) \geq \int r_T([a_0(x,w);x;w])
d\pi(x,w),$$
for all $T \geq T_0$.  As $\|R_T\|_F \leq \sqrt{\kappa}$ for all $T$, this sequence of integrals converges to zero.
Lemma~\ref{lem:a} now implies
\begin{flalign*}
  \int r_T([a_0(x,w)&;x;w]) d\pi(x,w)\\
  \geq\\\frac{1}{\d}\int
 |e_{\t_T}(a_0(x,w)) &- f_{\t_T}(x,w)|^2 d \pi(x,w).
\end{flalign*}
From that same lemma we see that the map $(x,w) \mapsto
 e_{\t_T}(a_0(x,w)) - f_{\t_T}(x,w)$ is real analytic.

 Let $L$ be the subspace of $\RR^{n_\t}$ defined by
 \begin{flalign*}
   L = \{\t \; | \;&\: 0 = e_\t(a_0(x,w)) - f_\t(x,w),
   \quad \forall \; (x,w) \in
   \X\times \W\}.
 \end{flalign*}
 Let $V$ be the quotient space $V = \RR^{n_\t}/L$, i.e. the set of
 equivalence classes where $\t \cong 0$ if 
 $e_\t(a_0(x,w))-f_\t(x,w)\equiv 0$.  The following functions
 are norms on $V$:
\begin{flalign*}
  \|\t\| =&\: \sqrt{\int|e_{\t}(a_0(x,w)) - f_{\t}(x,w)|^2d\pi(x,w)},\\
  |\t|_{\infty,U} =&\: \sup_{(x,w)\in U}\{|e_{\t}(a_0(x,w)) - f_{\t}(x,w)|\},
\end{flalign*}
where $U$ is any bounded open set.  Both functions are clearly
bounded, homogeneous, sub-additive and positive for all $\t$.  $\|\t\| =
|\t|_{\infty,U} = 0$ when $\t \cong 0$ so that the functions are
well-defined functions on $V$.  
As $U$ is open and $\pi$ is supported on an open set we see that $\t
\ncong 0$ implies both $\|\t\| > 0$ and $|\t|_{\infty,U} > 0$, 
by the Identity Theorem for analytic maps \cite{Fritzsche02}.  Now,
fix $U$ to be a bounded open set with $U \supset K$.  As all norms are
equivalent for finite dimensional spaces, we know there exists a
constant $c > 0$ such that $c\|\t\| \geq |\t|_{\infty,U}$.  That $\t_0
\cong 0$ implies $\|\t_T\| \cvgip 0$, so that 
$|\t_T|_{\infty,U} \cvgip 0$.  Lemma~\ref{lem:a} then yields
\begin{flalign*}
  \frac{2}{\sqrt{\d}}|\t_T|_{\infty,U}
 \geq \sup_{(x,w) \in K} \{ |e_{\t_T}^{-1}(f_{\t_T}(x,w))
  - a_0(x,w)|\},
\end{flalign*}
which establishes the claim.

\end{pf}

\bibliography{IEEEexample}{}
\bibliographystyle{plain}

\end{document}